\documentclass[a4paper,12pt]{article}
\usepackage{amsmath}
\usepackage{amsthm}
\usepackage{graphicx}
\usepackage{verbatim}
\usepackage{epsfig}
\usepackage{hyperref}
\usepackage{float}
\usepackage{color}
\usepackage{fullpage}
\usepackage{amsfonts}
\usepackage{amscd}
\usepackage{amssymb}
\usepackage{graphics}
\usepackage{amsmath}
\usepackage[english]{babel}
\usepackage{bbm}
\usepackage{yfonts}
\usepackage{mathrsfs}
\usepackage{color}
\newcounter{assumptions}

\newcommand{\BE}{{\mathbb{E}}}

\newcommand{\BP}{{\mathbb{P}}}

\newcommand{\BR}{{\mathbb{R}}}

\newcommand{\BZ}{{\mathbb{Z}}}

\newcommand{\FG}{{\mathfrak{G}}}

\newcommand{\CA}{{\mathcal{A}}}
\newcommand{\CB}{{\mathcal{B}}}
\newcommand{\CC}{{\mathcal{C}}}
\renewcommand{\CD}{{\mathcal{D}}}

\newcommand{\CF}{{\mathcal{F}}}

\newcommand{\CK}{{\mathcal{K}}}
\newcommand{\CL}{{\mathcal{L}}}

\newcommand{\CP}{{\mathcal{P}}}

\newcommand{\CW}{{\mathcal{W}}}

\newcommand{\diam}{{\rm diam}}
\newcommand{\Int}{{\rm Int}}

\newtheorem{theorem}{Theorem}[section]
\newtheorem{proposition}[theorem]{Proposition}
\newtheorem{lemma}[theorem]{Lemma}

\newtheorem{corollary}[theorem]{Corollary}

\newtheorem{example}[theorem]{Example}

\newtheorem{assumption}{Assumption}
\begin{document}
\numberwithin{equation}{section} \numberwithin{figure}{section}
\title{The existence phase transition for scale invariant 
Poisson random fractal models.}
\author{Erik I. Broman\footnote{Department of Mathematics, Chalmers University of Technology and Gothenburg University, Sweden. E-mail: broman@chalmers.se.
Supported by the Swedish Research Council}}
\maketitle
\begin{abstract}
In this paper we study the existence phase transition of scale invariant 
random fractal models. We determine the exact value of the critical point of
this phase transition for all models satisfying some weak assumptions.
In addition, we show that for a large subclass, the fractal model 
is in the empty phase at the critical point. This subclass of models 
includes the scale invariant Poisson
Boolean model and the Brownian loop soup.
In contrast to earlier results in the literature, we do not need to restrict 
our attention to random fractal models generated by open sets.
\end{abstract}
\section{Introduction}

The questions studied in this paper are very much related to the classical problem
of covering some fixed set by other randomly placed sets. The first 
paper dealing with this type of question was by 
Dvoretsky (\cite{Dvoretsky}) who studied the problem of 
covering the unit circle by independently and uniformly placing arcs 
of given lengths. The study of such covering problems later developed in several
different directions. One of the earliest papers dealing with questions 
similar to ours was by 
Shepp (\cite{Shepp}), who determined necessary and 
sufficient conditions for the real line to be covered by a Poisson 
process of random open intervals. The argument given in that paper relies
on the fact that the problem is posed in dimension one, and cannot 
readily be generalized to higher dimensions. Many years later, 
Bierm\'e and Estrade (\cite{BE}) and Broman, Jonasson and 
Tykesson (\cite{BJT}), studied the scale invariant Poisson Boolean model 
(or fractal ball model) in general dimensions. However, in those papers it was 
assumed that that the balls generating 
the fractal set were open, and in addition, 
the arguments given relied on the simple 
geometry of these balls.

The purpose of this paper is two-fold. Firstly, we want to find arguments
that do not rely on the sets generating the fractal being open. Secondly,
we want to give a general result that is not restricted to the case where
the generating sets have simple geometry (e.g. they are balls). 
Thereby, our results will cover cases such as the Brownian loop soup where
the outer boundary is closed and does not have a simple geometry (see 
also Example \ref{ex:BLS}).

Before stating exact results, we will begin by explaining the general setup.
Let $\FG$ be the set of subsets of $\BR^d$  
with diameter strictly smaller than 1.
Furthermore, let $\CF$ be some suitable $\sigma$-algebra
on $\FG.$ We consider infinite measures $\mu$ on $(\FG,\CF)$ which are 
semi scale invariant in the following sense. 
Assume that $\CD\in \CF$ is such that $\mu(\CD)<\infty,$ 
and let $0<s<\infty$ be such that 
\begin{equation} \label{eqn:CDs}
\CD_s:=\{G\subset \BR^d:G/s\in \CD\}
\end{equation}
only contains sets of diameter strictly smaller than one. 
Then, semi scale invariance simply means that $\mu(\CD_s)=\mu(\CD)$.
We shall also require $\mu$ to be translational invariant in that 
$\mu(x+\CD)=\mu(\CD)$ for every $\CD \in \CF$ and $x\in \BR^d.$ 
Here of course, 
$x+\CD=\{L\subset \BR^d:L=x+G \textrm{ for some } G\in \CD\}$. 

We then let $\lambda \mu$ (where $0<\lambda<\infty$ is a parameter) 
be the intensity measure of a Poisson process $\Phi_\lambda(\mu)$ on $\FG$.
Thus constructed, $\Phi_\lambda(\mu)$ is a (semi) scale and translation 
invariant random collection of bounded subsets of $\BR^d.$ This setup 
includes the scale invariant Poisson Boolean model (again, see \cite{BE} 
and \cite{BJT}) and the Brownian loop soup 
introduced by Lawler and Werner in \cite{LW}.

We do not give details of how to construct general triplets
$(\FG,\CF,\mu)$ in this paper, but instead we 
refer the reader to the book \cite{Molchanov} by Molchanov.
Our viewpoint is rather that given such a triplet, the results of 
this paper will apply. We also point out that there are plenty of 
such triplets including the examples mentioned above and the ones 
listed later in this paper.

\medskip

Throughout the rest of the paper, we shall make the following three 
natural assumptions on the measure $\mu.$ 
\begin{assumption} \label{assumption:reg1}
Let $A \subset \BR^d$ be bounded. Then, for any $\delta>0$ we have that 
\[
\mu(G\in \FG: G\cap A \neq \emptyset, \diam(G)\geq \delta)<\infty.
\]
\end{assumption}
\noindent
This assumption will make sure that the Poisson process never has 
infinitely many sets of diameter larger than some constant in any bounded region. 
Thus, Assumption 1 says that $\mu$ is locally finite. 

\medskip 

For our second assumption, let $\CL(\cdot)$ denote Lebesgue measure 
in $\BR^d$ and let $\partial G$ denote the boundary of the set $G$.
\begin{assumption} \label{assumption:reg2}
We have that $\mu(G\in \FG: \CL(\partial G)>0)=0$.
\end{assumption}
\noindent
This assumption makes sure that the boundary of the sets generating the fractal 
is not ''too large''.
\medskip

\begin{assumption} \label{assumption:reg3}
We have that $\mu(G\in \FG: \CL(G)>0)>0$.
\end{assumption}
\noindent
Let $\Int(G)$ denote the interior of a set $G.$
Note that under Assumption \ref{assumption:reg2} we have that 
$\CL(G \cup \partial G)=\CL(G)$, which means that if $\CL(G)=0$ then 
$\Int(G \cup \partial G)=\emptyset$ and so $G$ is a nowhere dense set. 
Therefore, if $\mu$ is concentrated on sets such that $\CL(G)=0$,
it follows from Baire's category theorem that $\CC(\lambda)\neq \emptyset$ 
for every $\lambda<\infty$. This explains the necessity of Assumption 
\ref{assumption:reg3} (see further the remark after Theorem 
\ref{thm:determinelambdae} below).

Note also that by Assumptions \ref{assumption:reg2} and \ref{assumption:reg3} 
we have that $\CL(\Int(G))\geq \CL(G)-\CL(\partial G)=\CL(G)>0$, 
and so $\mu$ puts positive mass on sets $G$ with non-empty interiors.

It will be convenient to split the measure $\mu$ into two parts. 
Firstly, we define the measure $\mu_p$ by letting 
$\mu_p(\CA)=\mu(\{G\in \CA:\CL(G)>0\})$ for any $\CA\in \CF$, 
and secondly we define 
$\mu_0$ by letting $\mu_0(\CA)=\mu(\{G\in \CA:\CL(G)=0\})$.
Obviously we then have that $\mu=\mu_p+\mu_0.$ 
Although we will consider the case where $\mu_0$ is allowed 
to be non-trivial, we point out that the case where $\mu=\mu_p$ is somewhat
easier to handle as it avoids some technical difficulties. 
See also the remark after the proof of Theorem \ref{thm:dieoutatcrit}.

\medskip

Our random fractal is then defined by letting
\begin{equation} \label{eqn:defCC}
\CC(\Phi_\lambda(\mu))
:=\BR^d \setminus \bigcup_{G\in \Phi_\lambda(\mu)} G,
\end{equation}
and we will write $\CC(\lambda)$ or simply $\CC.$ It is not hard to prove
(using Assumption \ref{assumption:reg3}) that for any $\lambda>0$ and 
any fixed $x\in \BR^d,$ we have that $\BP(x\in \CC)=0.$ 
We see that $\CC$ is a (semi) scale invariant random fractal, and we will 
be concerned with various properties of $\CC(\lambda)$ as $\lambda$ varies. 
It is useful to observe that by using a standard coupling argument, 
$\CC(\lambda)$ is decreasing in $\lambda.$

Random fractal models exhibit several phase transitions 
(see for instance Dekking and Meester, \cite{DM}). The two 
most studied are the {\em existence} and the
{\em connectivity} phase transitions, as we now explain.
Define
\[
\lambda_e:=\inf\{\lambda>0:\BP(\CC(\lambda)=\emptyset)=1\}.
\]
Thus for $\lambda>\lambda_e,$ $\CC(\lambda)$ is almost surely empty,
and we say that it is in the empty phase. Clearly, this is the same
as saying that $\BR^d$ is covered by $\bigcup_{G\in \Phi_\lambda(\mu)} G.$ 
If on the other hand $\lambda<\lambda_e,$ then 
$\BP(\CC(\lambda)\cap [0,1]^d\neq \emptyset)>0.$ Thus, by tiling $\BR^d$
with translates of $[0,1]^d$ and noting that 
$\CC(\lambda) \cap B_1$ and $\CC(\lambda) \cap B_2$ are independent 
whenever $B_1,B_2 \subset \BR^d$ are such that $d(B_1,B_2)\geq 1,$ we conclude
that $\BP(\CC(\lambda)\neq \emptyset)=1$. 
Hence, $\lambda_e$ is the critical point of the existence phase 
transition. 
Furthermore, we can define
\[
\lambda_c:=\sup\{\lambda>0:\BP(\CC(\lambda)
\textrm{ contains connected components larger than one point})=1\}.
\]
This means that if $\lambda>\lambda_c,$ then $\CC(\lambda)$ is almost surely
totally disconnected. However, if $\lambda<\lambda_c,$ then $\CC(\lambda)$
will almost surely contain connected components (by the same tiling 
argument as before).

It is natural and interesting to ask what happens at the critical points
of these phase transitions. In \cite{BC} it was proven that
\[
\BP(\CC(\lambda_c)
\textrm{ contains connected components larger than one point})=1,
\]
for any random fractal model satisfying Assumptions \ref{assumption:reg1}
and \ref{assumption:reg3}. Therefore,
{\em at} $\lambda_c$ the corresponding fractal is in the connected phase.
In that respect, this phase transition is very well understood.

As explained, the existence phase transition is closely related to 
classical problems of random coverings, and we can now explain the results of
\cite{BE} and \cite{BJT} in terms of the terminology just introduced. 
In \cite{BE}, the exact value of $\lambda_e$ was determined for the 
fractal ball model.
In fact, their result was stronger than this, but it did not 
cover what happened at the critical value, i.e. whether 
$\BP(\CC(\lambda_e)=\emptyset)=0$ or $1$. This case was determined in \cite{BJT} 
where it was proven that $\BP(\CC(\lambda_e)=\emptyset)=1$ for the particular
case when $\mu$ is supported on open balls. In that paper, the corresponding
result was also proven when  $\mu$ is supported on open boxes.

Let 
\begin{equation} \label{eqn:defAl}
\CA_l:=\{G \in \FG: o\in G, 2^{-l} \leq \diam(G)< 2^{-l+1}\},
\end{equation}
where $o$ denotes the origin in $\BR^d$. It follows by scale 
invariance that $\mu(\CA_l)=\mu(\CA_{l+1})$ for every $l\geq 1.$ \\
\noindent
{\bf Remark:}
In light of the definition of $\CD_s$ in \eqref{eqn:CDs}
above, it might be more natural to denote the left hand side of \eqref{eqn:defAl} 
by $\CA_{2^{-l}}.$ However, it would quickly be cumbersome to use this notation.
Similar remarks apply throughout the paper.

\medskip

Our first main result is the following. 
\begin{theorem} \label{thm:determinelambdae}
For any $\mu$ satisfying Assumptions \ref{assumption:reg1},
\ref{assumption:reg2} and \ref{assumption:reg3} we have that 
\[
\lambda_e=\frac{d \log 2}{ \mu(\CA_1)}.
\]
\end{theorem}
\noindent
{\bf Remarks:} 
The use of $\CA_1$ is somewhat arbitrary. Indeed, by scaling invariance
we have that for any $k \geq 1,$
\[
\mu\left(G \in \FG: o\in G, 2^{-1}\leq \diam(G) < 1\right)
=k\mu\left(G \in \FG: o\in G, 2^{-1/k} \leq \diam(G) < 1\right).
\]
%\begin{eqnarray*}
%\lefteqn{\mu\left(G \in \FG: o\in G, 2^{-1}\leq \diam(G) < 1\right)}\\
%& &=k\mu\left(G \in \FG: o\in G, 2^{-1/k} \leq \diam(G) < 1\right).
%\end{eqnarray*}
As we shall see in Section \ref{sec:proofcrit}, there is a 
canonical way of rewriting the expression for $\lambda_e$. Informally,
this alternative expression involves the (quasi-)expected volume of 
sets of diameter exactly one.
However, as this discussion is too long to fit in the introduction 
we defer it to Section \ref{sec:proofcrit}
where the alternative expression appears in \eqref{eqn:lambdae_alt}.

From alternative expressions for $\mu(\CA_1)$ derived later in the paper 
(i.e. \ref{eqn:muA1alt}), it is easy to see that if Assumption 
\ref{assumption:reg3} is not satisfied, (i.e. if 
$\mu(G\in \FG:\CL(G)>0)=0$), then it follows that $\mu(\CA_1)=0.$ 
According to the discussion after the statement of 
Assumption \ref{assumption:reg3}, we will in this situation have
that $\lambda_e=\infty.$ If we then just interpret $d\log 2/0$ as 
being $\infty$, it is possible to drop Assumption \ref{assumption:reg3} 
from the statement of Theorem \ref{thm:determinelambdae}.

It is natural to ask whether Theorem \ref{thm:determinelambdae} still holds
when Assumption \ref{assumption:reg2} is not satisfied. We do not have any example
showing that this is not the case.

\bigskip

We now turn to the behaviour of $\CC(\lambda)$ at the critical point $\lambda_e.$
For any set $G,$ let $[\partial G]^r=\{x\in G: d(x, \partial G ) < r\} \subset G$,
so that $[\partial G]^r$ is the inner $r$-neighbourhood of $\partial G$. Our
second main result is the following.
\begin{theorem} \label{thm:dieoutatcrit}
If $\mu$ satisfies Assumptions \ref{assumption:reg1}, \ref{assumption:reg2},
\ref{assumption:reg3}
and the additional condition that
\begin{equation} \label{eqn:extracond}
\int_{\CA_1} \int_0^1 \frac{1}{r} \frac{\CL([\partial G]^r)}{\CL(G)} 
dr d\mu_p(G)<\infty,
\end{equation}
%\sum_{l=1}^\infty \int_{\CA_1}
%\frac{\CL(G\setminus S(G,2^{-l}))}{\CL(G)} d\mu(G)<\infty,
then we have that $\CC(\lambda_e)=\emptyset$ almost surely.
\end{theorem}
\noindent
{\bf Remark:}
The condition \eqref{eqn:extracond} in Theorem \ref{thm:dieoutatcrit} 
might look like a technical 
assumption made purely for convenience. However, the class of models that 
satisfy \eqref{eqn:extracond} is very large and includes the fractal 
ball model and the fractal generated by the outer boundaries of the 
Brownian loop soup. As we will see, our examples of models that do not 
satisfy Assumption \ref{assumption:reg2} or condition \eqref{eqn:extracond}
will be somewhat contrived. This of course does not rule out the possibility 
of a naturally occurring example which does not satisfy our assumptions 
and conditions. Note also the occurrence of $\mu_p$ in \eqref{eqn:extracond}
which guarantees that $\CL(G)>0$ in the integrand.

The above mentioned examples will briefly be 
analysed in Section \ref{sec:examples} (along with a third example) where we 
discuss what our main theorems imply for them.
It is appropriate to mention the paper by Nacu and Werner (\cite{NW}), 
which deals with Hausdorff dimensions of the so-called Gasket of 
two-dimensional loop soups. In that paper a condition appears under 
the name ''thin'', which is a slightly 
stronger version of condition \eqref{eqn:extracond}. 
We give the exact statement of this condition in \eqref{eqn:thin}.
We note that condition \eqref{eqn:extracond}, Assumption \ref{assumption:reg2}
and the condition of being thin all deal with the boundaries
$\partial G$ of the sets used to generate the random fractal.
It is therefore natural to ask how Assumption \ref{assumption:reg2} 
and conditions \eqref{eqn:extracond}, \eqref{eqn:thin} relate 
to each other. The answer to this is given by
the following proposition.

\begin{proposition} \label{prop:condassrel}
 \

\begin{itemize}

\item[(i)] If the measure $\mu$ is thin, then Assumption %s \ref{assumption:reg1}, 
\ref{assumption:reg2} and condition \eqref{eqn:extracond} are both satisfied. 

\item[(ii)] If $\mu$ is supported on compact sets, then
condition \eqref{eqn:extracond} implies Assumption \ref{assumption:reg2},
while for general measures $\mu$ this is not the case.

\item[(iii)] Assumption \ref{assumption:reg2} does not imply 
condition \eqref{eqn:extracond}.

\end{itemize}
\end{proposition}
\noindent
{\bf Remark:} Assumption \ref{assumption:reg1} is a fundamentally different type 
of condition than the ones we deal with in Proposition \ref{prop:condassrel}, 
and it is therefore not so natural to try to relate it to these
other conditions. However, we add a short discussion on this subject 
after the proof of Proposition \ref{prop:condassrel}.

We also note that the proof of Proposition \ref{prop:condassrel} provides 
a concrete example of a measure $\mu$ not covered by our results. It will 
be a measure that does not satisfy Assumption \ref{assumption:reg2}.

\bigskip

We will mention one more result. 
Assume momentarily that $\mu$ is supported on open sets, and that $\tilde{\mu}$
is the measure induced by taking the closure of these sets. We then have the 
following corollary of Theorems \ref{thm:determinelambdae} and 
\ref{thm:dieoutatcrit}.
\begin{corollary} \label{corr:closedopensame}
If $\mu$ satisfies Assumptions \ref{assumption:reg1}, \ref{assumption:reg2} 
and \ref{assumption:reg3}, then 
the critical values for $\mu$ and $\tilde{\mu}$ are the same, i.e. 
$\lambda_e(\mu)=\lambda_e(\tilde{\mu}).$ 
\end{corollary}
\noindent
In particular, Corollary \ref{corr:closedopensame} implies that the two 
fractal ball models generated by open balls
and closed balls respectively, have the same critical value. Furthermore,
it follows from Theorem \ref{thm:dieoutatcrit} (see Example \ref{ex:ball}) 
that $\CC(\lambda_e)=\emptyset$ for both of these variants.

\medskip

The rest of the paper is organized as follows. In Section \ref{sec:defmodels}
we establish necessary notation and a preliminary result.
In Section \ref{sec:proffthm1} we will prove Theorem 
\ref{thm:determinelambdae} while in Section \ref{sec:proofcrit} we 
prove Theorem \ref{thm:dieoutatcrit}. In addition, Section \ref{sec:proofcrit} 
also contains the proof of Proposition \ref{prop:condassrel} and
the alternative expression for $\lambda_e$ mentioned above. This alternative 
expression will be used to prove Corollary \ref{corr:closedopensame}.
Finally, Section \ref{sec:examples} contains the above mentioned examples.

\section{Models and preliminary results} \label{sec:defmodels}

As explained in the introduction, the measure $\mu$ is a measure on 
\[
\FG:=\{G\subset \BR^d: \diam(G) < 1\},
\]
and $\Phi_\lambda=\Phi_\lambda(\mu)$ is a Poisson process on 
$\FG$ using $\lambda \mu$ as its intensity measure.
Furthermore, we let 
\[
\FG_n:=\{G\subset \BR^d: 2^{-n} \leq \diam(G) < 1\},
\]
and
\[
\Phi_{\lambda,n}:=\{G\in \Phi_\lambda (\mu): G \in \FG_n\},
\]
so that $\Phi_{\lambda,n}$ is a subset of the process $\Phi_\lambda$ containing 
sets of diameter
larger than or equal to $2^{-n}.$ Then, define
\[
\CC_n:=\BR^d \setminus \bigcup_{G \in \Phi_{\lambda,n}} G,
\] 
and note that $\CC_n \supset \CC_{n+1}$ for every $n.$ Note also that by 
\eqref{eqn:defCC} $\CC_n \downarrow \CC.$
%\downarrow \CC.$ 
For $m>n,$ let
\[
\CC_{m}^n:=\BR^d \setminus 
\bigcup_{G \in \Phi_{\lambda,m}\setminus \Phi_{\lambda,n}} G,
\]
so that $\CC_m^n \cap \CC_n=\CC_m$, and $\CC_m^n, \CC_n$ 
%have no $K \in \Phi_{\lambda}(\mu)$ in common and 
are independent. 
Next, let 
\[
\FG^o:=\{G\in \FG: o\in G\},
\]
and
\[
\FG^o_n:=\{G\in \FG_n: o\in G\},
\]
so that $\FG_n^o$ consists of the sets that have diameter larger than 
$2^{-n}$ while containing the origin. 
We observe from \eqref{eqn:defAl} that 
\[
\FG^o_n=%\{G\in \FG: \diam(G)=1\}
\bigcup_{l=1}^n \CA_l,
\]
and by scaling invariance we have that 
\[
\mu(\FG_n^o)=\sum_{l=1}^n \mu(\CA_l)=n \mu(\CA_1).
\]

Next, we let
\[
\CB_n:=\left\{x+[0,2^{-n}]^d:
x \in (2^{-n}\BZ^d)\cap [0,1-2^{-n}]^d\right\},
\]
so that $\CB_n$ consists of $2^{dn}$ closed boxes with non-overlapping
interiors and $\bigcup_{B \in \CB_n}B=[0,1]^d$. 
We shall refer to a box $B\in \CB_n$ as a level $n$ box.
Furthermore, for any $A \subset \CB_n$ we define 
$|A|$ to be the cardinality of $\{B\in \CB_n:B \in A\}$.
Then, let 
\[
M_n:=\{B\in \CB_n:
\not \exists \ G \in \Phi_{\lambda,n} \textrm{ such that } B \subset G \}.
\]
Of course, we could have used $\Phi_\lambda$ in the definition without 
anything being changed.
Thus, $M_n$ is the set of level $n$ boxes which are not singly covered by a 
set in the Poisson process $\Phi_{\lambda,n}.$ Furthermore, let
\[
m_n:=\{B\in \CB_n:
\not \exists \ G \in \Phi_{\lambda,n} \textrm{ such that } B\cap G\neq \emptyset \}.
\]
Hence, $m_n$ is the set of level $n$ boxes that are completely untouched
by the sets of $\Phi_{\lambda,n}.$ It is immediate that $m_n \subset M_n$ 
for every $n$. Let $D_n=D_n(\CC_n)$ be a minimal (with respect to the 
number of boxes) collection of boxes in $\CB_n$ such that 
\[
\CC_n\cap [0,1]^d \subset \bigcup_{B\in D_n} B.
\]
Since a point $x\in \CC_n$ in the intersection of two boxes $B_1$ and 
$B_2$ can be covered by either one of them,
it follows that $D_n$ is not necessarily unique. 
If there is more than one way of choosing such a set $D_n,$ 
we use some predetermined rule to pick one.
Finally, we let $L_n:=|D_n|$.
%, i.e. the number of level $n$ boxes needed to cover 
%$\CC_n \cap [0,1]^d.$ %|\{B\in D_n\}|$.

A version of the following easy lemma appears in \cite{BJT}. However, 
since the proof is short, we include it here for completeness. 

\begin{lemma} \label{lemma:survMninf}
For any $\mu$ satisfying Assumption \ref{assumption:reg3}, 
if $\lambda>0$ is such that 
$\BP(\CC(\lambda) \neq \emptyset)=1$, then we must have that 
\[
\BP(\lim_{n \to \infty}|M_n|=\infty)>0.
\] 
\end{lemma}
\noindent
{\bf Proof.}
First, any $B\in D_n$ contains a point which is in $\CC_n(\lambda)$ and 
therefore it cannot be that $B$ is covered by a single 
set $G$ in the Poisson process $\Phi_{\lambda,n}$. Therefore, $B\in M_n$ and so  
$L_n \leq |M_n|.$

Secondly, 
%let
%\[
%E_{n}:=\bigcup_{B \in D_{n}} B ,
%\]
%and 
observe that by definition of $D_{n}$ we have that
\[
(\CC_{n}\cap[0,1]^d)\setminus \bigcup_{B \in D_{n}} B=\emptyset.
\]
%By assumption, every $K\in \Phi_\lambda(\mu)$ has non-empty interior, and so
Furthermore, we have that by Assumption \ref{assumption:reg3} there 
exists some $\alpha=\alpha(\lambda)\in(0,1),$ such that 
\[
\BP(\CC_1\cap [0,1]^d=\emptyset)=\alpha.
\]
By using the FKG inequality for Poisson processes together
with the scaling invariance of the models and the fact that $D_n$ is measurable 
with respect to the $\sigma$-algebra generated by $\Phi_{\lambda,n}$ ,
we conclude that
\begin{eqnarray} \label{eqn:alphaL}
\lefteqn{\BP(\CC_{n+1}\cap [0,1]^d=\emptyset \ | \ \Phi_{\lambda,n})}\\
& & \geq \BP\left(\bigcap_{B \in D_{n}} \{\CC_{n+1}^{n} \cap B=\emptyset\} 
\ \Bigg{|} \ \Phi_{\lambda,n}\right)
\geq \prod_{B\in D_{n}}
\BP(\CC_{n+1}^{n}\cap B=\emptyset)
= \alpha^{L_{n}}>0. \nonumber
\end{eqnarray}
Therefore, if there exists $L<\infty$ such that
$L_{n}\leq L$ for infinitely many $n,$ 
it follows by standard arguments that \eqref{eqn:alphaL} implies that
$\CC\cap [0,1]^d=\emptyset$  almost surely. By countability, we conclude
that if $\BP(\CC \neq \emptyset)=1$ we must have that 
$\BP(\lim_{n \to \infty} L_n=\infty)>0.$  Since $L_n \leq |M_n|$ the 
statement follows.
\fbox{}\\
% SEE FOR INSTANCE https://math.stackexchange.com/questions/312151/generalized-second-borel-cantelli-lemma FOR SORT OF THIS THING.

\section{Proof of Theorem \ref{thm:determinelambdae}} \label{sec:proffthm1}
In this section, we will prove Theorem \ref{thm:determinelambdae}. We will do
this by showing the two directions separately.

\subsection{Lower bound}
This subsection is devoted to proving that 
$\lambda_e \geq \frac{d \log 2}{\mu(\CA_1)},$ and this is 
done in two steps. Firstly, Lemma \ref{lemma:small_lambda} proves
that for any $\lambda<\frac{d \log 2}{\mu(\CA_1)}$ we must have that 
$\BE[|m_n|] \to \infty.$ Then, in Theorem \ref{thm:lbound} we use this 
to embed a supercritical branching process of boxes belonging to 
$(m_n)_{n=1}^\infty$ 
in such a way that
$\CC(\lambda)\cap [0,1]^d \neq \emptyset$ whenever the branching 
process survives. In fact, with minimal extra effort, we will obtain 
the stronger statement that for every $\epsilon>0$ we can embed a 
supercritical branching process with growth rate 
(see \cite{LP} Chapter 1) at least
$2^{\left(d-\lambda \frac{\mu(\CA_1)+\epsilon}{\log 2}\right)}$.

\begin{lemma} \label{lemma:small_lambda}
For any $\lambda,\epsilon>0$ and fractal model satisfying Assumptions 
\ref{assumption:reg1},\ref{assumption:reg2} and \ref{assumption:reg3}
there exists $N=N(\epsilon)<\infty$ such that 
\[
\BE[|m_n|]\geq 2^{n\left(d-\lambda \frac{\mu(\CA_1)+\epsilon}{\log 2}\right)},
\]
for every $n\geq N.$
\end{lemma}
\noindent
{\bf Proof.}
Let $D(x,r)$ denote an open ball centered at $x$ with radius $r.$
Obviously, the box $[0,2^{-n}]^d$ can be inscribed in an open ball of 
radius $d2^{-n}.$ Therefore, by translation invariance, 
we have that for any $B\in \CB_n,$ 
\begin{eqnarray} \label{eqn:ref1}
\lefteqn{\BP(B\in m_n)\geq 
\BP(\not \exists \ G\in \Phi_{\lambda,n}: D(o,d2^{-n})\cap G \neq \emptyset)}\\
& & =\exp\left(-\lambda \mu(\{G\in \FG_n:\ D(o,d2^{-n}) \cap G\neq \emptyset\})\right). \nonumber 
\end{eqnarray}
Let 
\[
E(G,r):=\{x:d(x,G)< r\}
\]
so that $E(G,r)$ is an enlargement of $G.$ If $o\in G$ then trivially 
$D(o,r)\cap G \neq \emptyset$ and $o\in E(G,r).$ Furthermore, if $o \in G^c,$
then since $D(o,r)$ and $E(G,r)$
are both open we see that $D(o,r)\cap G \neq \emptyset$ implies that 
$d(o,\partial G)<r$ and so $o\in E(G,r).$ Furthermore, the reverse statements
are also true. Therefore,
\begin{eqnarray}\label{eqn:muKn1}
\lefteqn{\mu(\{G\in \FG_n: D(o,d2^{-n})\cap G \neq \emptyset \})}\\
& & =\mu(\{G\in\FG_n: o\in E(G,d2^{-n})\})
=\int_{\FG_n} I(o\in E(G,d2^{-n})) d\mu(G). \nonumber
\end{eqnarray}
Note that for any fixed point $x\in \BR^d$ we must have that
\begin{eqnarray*}
\lefteqn{-x+\{G\in \FG_n:o\in E(G,d2^{-n}),x\in G\}}\\
& & =\{L\in \FG_n:L=-x+G, o\in E(G,d2^{-n}),x\in G\} \\
& & =\{G\in \FG_n:o\in E(G+x,d2^{-n}),x\in (G+x)\}.
%& & =\{L\in \FG_n:-x \in E(L,d2^{-n}),o\in L\} 
\end{eqnarray*}
Recall that $\CL(\cdot)$ denotes Lebesgue measure in $\BR^d$
and recall also the notation $\mu=\mu_p+\mu_0$. 
At this point we will need to consider $\mu_p$ and $\mu_0$
separately. We will start with $\mu_p$ and 
note that by translation invariance we have that
\begin{eqnarray} \label{eqn:Eequal}
\lefteqn{\int_{\FG_n}I(o\in E(G,d2^{-n}))d\mu_p(G)}\\
& & =\int_{\FG_n}I(o\in E(G,d2^{-n}))\frac{1}{\CL(G)}
\int_{\BR^d}I(x\in G)dx d\mu_p(G)\nonumber \\
& & =\int_{\BR^d}\int_{\FG_n}\frac{1}{\CL(G)}
I(o\in E(G,d2^{-n}))I(x\in G) d\mu_p(G) dx\nonumber \\
& & =\int_{\BR^d}\int_{\FG_n}\frac{1}{\CL(G+x)}
I(o\in E(G+x,d2^{-n})I(x\in (G+x)) d\mu_p(G) dx\nonumber\\
& & =\int_{\BR^d}\int_{\FG_n}\frac{1}{\CL(G)}
I(o\in E(G,d2^{-n})+x)I(o\in G) d\mu_p(G) dx \nonumber\\
& & =\int_{\BR^d}\int_{\FG_n}\frac{1}{\CL(G)}
I(-x\in E(G,d2^{-n}))I(o\in G) d\mu_p(G) dx \nonumber\\
%& & =\int_{\FG_n} \frac{1}{\CL(G)}
%I(o\in G)\int_{\BR^d}I(-x\in E(G,d2^{-n})) dx d\mu(G) \\
& & =\int_{\FG_n} \frac{\CL(E(G,d2^{-n}))}{\CL(G)}
I(o\in G)d\mu_p(G) \nonumber \\
& & 
=\int_{\FG^o_n} \frac{\CL(E(G,d2^{-n}))}{\CL(G)}d\mu_p(G)
=\sum_{l=1}^n\int_{\CA_l} \frac{\CL(E(G,d2^{-n}))}{\CL(G)}d\mu_p(G). \nonumber
\end{eqnarray}
We make the observation that since by \eqref{eqn:muKn1} and \eqref{eqn:Eequal},
\begin{equation} \label{eqn:Efinite}
\sum_{l=1}^n\int_{\CA_l} \frac{\CL(E(G,d2^{-n}))}{\CL(G)}d\mu_p(G)
\leq \mu(\{G\in \FG_n:\ D(o,d2^{-n}) \cap G\neq \emptyset\}),
\end{equation}
then the left hand side must be finite by Assumption \ref{assumption:reg1}. 
Then, let  
\[
a_{l,n}:=\int_{\CA_l} \frac{\CL(E(G,d2^{-n}))}{\CL(G)}d\mu_p(G),
\]
and observe that by elementary properties of Lebesgue measure we have that 
\[
2^d\CL(E(G,d2^{-n}))
=\CL(2E(G,d2^{-n}))
=\CL(E(2G,d2^{-n+1})).
\]
Therefore, 
\begin{eqnarray} \label{eqn:alnrel}
\lefteqn{a_{l,n}=\int_{\CA_l} \frac{\CL(E(G,d2^{-n}))}{\CL(G)}d\mu_p(G)}\\
& & =\int_{\CA_l} \frac{\CL(E(2G,d2^{-n+1}))}{\CL(2G)}d\mu_p(G)
=\int_{\CA_{l-1}} \frac{\CL(E(G,d2^{-n+1}))}{\CL(G)}d\mu_p(G)=a_{l-1,n-1}. \nonumber
\end{eqnarray}
Here, the third equality follows since the measure $\mu$ is scale invariant,
and the fact that for any $G\in \CA_l$ we have that 
$2G\in \CA_{l-1}$.

By Assumption \ref{assumption:reg2} we have that
$\mu(G\in \FG: \CL(\partial G )>0)=0,$
and so 
\[
a_{1,n}=\int_{\CA_1} \frac{\CL(E(G,d2^{-n}))}{\CL(G)}
I(\CL(\partial G)=0)d\mu_p(G).
\]
Furthermore, we have that 
\[
\bigcap_{n=1}^\infty E(G,d2^{-n})=G \cup \partial G,
\]
and so by the monotone convergence theorem we conclude that 
\begin{eqnarray} \label{eqn:b1nconverge}
\lefteqn{\lim_{n \to \infty} a_{1,n}
=\int_{\CA_1} \lim_{n \to \infty}
\frac{\CL(E(G,d2^{-n}))}{\CL(G)}I(\CL(\partial G)=0)d\mu_p(G)}\\
& & =\int_{\CA_1} 
\frac{\CL(G \cup \partial G)}{\CL(G \cup \partial G)}
I(\CL(\partial G)=0)d\mu_p(G)
=\mu_p(\CA_1). \nonumber
\end{eqnarray}
By \eqref{eqn:Efinite}, and the elementary fact that $a_{1,n}\geq a_{1,n+1}$
for any $n$, there exists some 
$C_\mu<\infty$ such that $a_{1,n}\leq C_\mu$ for every $n \geq 1.$ 
For any fixed $\epsilon>0$ 
it follows from \eqref{eqn:b1nconverge} that there exists an 
$N_1<\infty$ such that 
$a_{1,n}\leq \mu(\CA_1)+\epsilon/3$ for every $n \geq N_1.$ Furthermore, 
there exists $N \geq N_1$ such that 
\[
(N-N_1)(\mu_p(\CA_1)+\epsilon/3)+C_\mu N_1 \leq 
N(\mu_p(\CA_1)+\epsilon/2).
\]
By \eqref{eqn:alnrel} we have that $a_{l,n}=a_{l-1,n-1}=\cdots=a_{1,n-l+1}$, and
so we see that for any $n\geq N,$
\begin{eqnarray} \label{eqn:alnsum}
\lefteqn{\sum_{l=1}^n a_{l,n}=\sum_{l=1}^n a_{1,n-l+1}
=\sum_{l=1}^n a_{1,l}}\\
& & \leq \sum_{l=1}^{N_1} a_{1,l} +\sum_{l=N_1+1}^n a_{1,l}
\leq \sum_{l=1}^{N_1} C_\mu +\sum_{l=N_1+1}^n (\mu_p(\CA_1)+\epsilon/3) 
\nonumber \\
& &=C_\mu N_1+(n-N_1)(\mu_p(\CA_1)+\epsilon/3)
\leq n(\mu_p(\CA_1)+\epsilon/2). \nonumber
\end{eqnarray}
Combining \eqref{eqn:Eequal}, \eqref{eqn:alnrel} and \eqref{eqn:alnsum}
we obtain 
\begin{equation} \label{eqn:mupineq}
\int_{\FG_n}I(o\in E(G,d2^{-n}))d\mu_p(G)
\leq n(\mu_p(\CA_1)+\epsilon/2),
\end{equation}
for every $n\geq N.$ 

We are now ready to turn our attention to $\mu_0.$ However, as this case
is similar to our previous case we will be brief.
We have that (with the interpretation that $\FG_0=\emptyset$)
\begin{eqnarray*} 
\lefteqn{\int_{\FG_n}I(o\in E(G,d2^{-n}))d\mu_0(G)}\\
& & =\sum_{l=1}^n \int_{\FG_l\setminus \FG_{l-1}}I(o\in E(G,d2^{-n}))d\mu_0(G)
=\sum_{l=1}^n \int_{\FG_1}I(o\in E(G,d2^{-n+l-1}))d\mu_0(G) \\
& & =\sum_{l=1}^n \int_{\FG_1}I(o\in E(G,d2^{-n+l-1}))
\frac{1}{\CL(E(G,1))}\int_{\BR^d}I(x\in E(G,1))dx d\mu_0(G) \\
& & =\sum_{l=1}^n \int_{\FG_1}
\frac{\CL(E(G,d2^{-n+l-1}))}{\CL(E(G,1))}I(o\in E(G,1))d\mu_0(G)\\
& & =\sum_{l=1}^n \int_{\FG_1}
\frac{\CL(E(G,d2^{-l}))}{\CL(E(G,1))}I(o\in E(G,1))d\mu_0(G).
\end{eqnarray*}
Here we used the scale invariance in the second equality
and then a calculation similar to \eqref{eqn:Eequal}.
Next, note that by Assumption \ref{assumption:reg1} and the monotone 
convergence theorem we have that 
\begin{eqnarray*} 
\lefteqn{\lim_{l \to \infty} \int_{\FG_1}
\frac{\CL(E(G,d2^{-l}))}{\CL(E(G,1))}I(o\in E(G,1))d\mu_0(G)}\\
& & =\int_{\FG_1}
\frac{\CL(G\cup \partial G)}{\CL(E(G,1))}I(o\in E(G,1))d\mu_0(G)=0,
\end{eqnarray*}
since $\CL(G\cup \partial G)=0$ on the support of $\mu_0.$
Therefore, by taking  $N$ perhaps even larger than before, we obtain that
\begin{equation} \label{eqn:mu0ineq}
\int_{\FG_n}I(o\in E(G,d2^{-n}))d\mu_0(G)
\leq n\epsilon/2,
\end{equation}
for every $n \geq N.$ We note that (as in the above calculations)
\[
\mu_0(\CA_1)=\int_{\FG_1} I(o \in G) d\mu_0(G)
\leq \lim_{l\to \infty} \int_{\FG_1} I(o \in E(G,d 2^{-l})) d\mu_0(G)=0.
\]
Hence, $\mu(\CA_1)=\mu_p(\CA_1)+\mu_0(\CA_1)
=\mu_p(\CA_1)$ and so by combining \eqref{eqn:ref1}, \eqref{eqn:muKn1},
\eqref{eqn:mupineq} and \eqref{eqn:mu0ineq} we then see that for 
every $n \geq N,$
\begin{eqnarray*}
\lefteqn{\BE[|m_n|]
\geq2^{dn}\exp\left(-\lambda\int_{\FG_n}I(o\in E(G,d2^{-n}))d\mu(G)\right)}\\
& & =2^{dn}\exp\left(-\lambda \int_{\FG_n}I(o\in E(G,d2^{-n}))d\mu_p(G) 
-\lambda \int_{\FG_n}I(o\in E(G,d2^{-n}))d\mu_0(G) \right)\\
& & \geq 2^{dn}\exp\left(-\lambda n(\mu_p(\CA_1)+\epsilon)\right)
=2^{n\left(d-\lambda \frac{\mu(\CA_1)+\epsilon}{\log 2}\right)}
\end{eqnarray*}
as desired.
\fbox{}\\
\noindent
{\bf Remark:} Note the we do not really use Assumption \ref{assumption:reg3}
in the above proof. The only purpose of this assumption is to guarantee that 
$\mu(\CA_1)=\mu_p(\CA_1) > 0.$

\medskip

Let 
\[
\FG_{lN}^{(l-1)N}:=\{G \in \FG:2^{-lN}\leq \diam(G)<2^{-(l-1)N}\},
\]
and 
\[
\Phi_{lN}^{(l-1)N}:=\{G \in \Phi_\lambda(\mu): G \in \FG_{lN}^{(l-1)N} \}.
\]

\begin{theorem} \label{thm:lbound}
For any fractal model satisfying Assumptions \ref{assumption:reg1},
\ref{assumption:reg2} and \ref{assumption:reg3} we have that  
\[
\lambda_e \geq \frac{d \log 2}{\mu(\CA_1)}.
\]
\end{theorem}
\noindent
{\bf Proof.}
The strategy is to embed the family tree of a supercritical branching process 
within the sequence $(m_n)_{n \geq 1}$. The boxes of this branching process 
will then be used to prove that $\CC(\lambda)\cap [0,1]^d \neq \emptyset$
whenever the branching process survives.

Fix $\lambda<\frac{d \log 2}{\mu(\CA_1)}$ and $\epsilon>0$ such that 
$\alpha:=d-\lambda\frac{\mu(\CA_1)+\epsilon}{\log 2}>0$.
The first step is to pick $N$ so that $\BE[|m_{N}|]\geq 3^d 2^{\alpha N}$ 
which we can do by 
Lemma \ref{lemma:small_lambda}.
Given the set $m_N \subset \CB_N$, we then let $W_1=W_1(m_N)$ be a maximal subset
such that $d(B_i,B_j)\geq 2^{-N}$ for every $B_i,B_j\in m_N$ with $i \neq j.$
For definitiveness, we can 
assume that $W_1$ is picked by some predetermined rule. Assume also that 
$\{B^1_1,B^1_2,\ldots,B^1_{Z_1}\}$ is an enumeration of the boxes $B\in W_1$ 
where of course, $Z_1:=|W_1|.$  We observe that $Z_1\geq m_N/3^d$. This follows  
since any two level $N$ boxes $B,B'\in \CB_N$ which have no points in common 
(i.e. they are not 
neighbours), must be at distance at least $2^{-N}$ from each other. 
We see that 
\[
\BE[Z_1]\geq \frac{\BE[|m_N|]}{3^d}\geq 2^{\alpha N},
\]
by our assumption on $N.$ Obviously, $Z_1 \leq |m_N|$ and $Z_1$ is the 
first generation of our branching process, while the set of boxes in $W_1$
is the first generation of the corresponding family tree. Finally, let 
\[
\CW_1:=\bigcup_{B \in W_1} B,
\]
and observe that $\CW_1 \subset \CC_N.$

We proceed with the construction of the second generation.
For any fixed $B^1_i\in\{B^1_1,\ldots,B^1_{Z_1}\}$, let 
\[
\CB_{2N}(B^1_i)=\{B \in \CB_{2N}: B \subset B^1_i\}
\]
and let
\[
m_{2N,N}^i:=\{B\in \CB_{2N}(B^1_i):
\not \exists \ G \in \Phi_{2N}^N \textrm{ such that } B\cap G\neq \emptyset \}.
\]
We note that a box $B\in m_{2N,N}^i$ is untouched by any set (in $\Phi_\lambda$) 
of diameter 
between $2^{-2N}$ and $2^{-N}$. In addition, since $B\subset \CW_1$ it must 
in fact also be untouched by any set (in $\Phi_\lambda$) 
of diameter larger than $2^{-N}$ and so $B \subset \CC_{2N}$.

We now make two key observations. 
Firstly, since $d(B^1_i,B^1_j)\geq 2^{-N}$ whenever $i\neq j,$ we have that 
\[
\{G\in \FG_{2N}^N: G\cap B^1_i \neq \emptyset\}
\cap \{G\in \FG_{2N}^N: G\cap B^1_j \neq \emptyset\}
=\emptyset
\]
for any $i \neq j.$ This of course follows since for any 
$G \in \FG_{2N}^N$ we have that $\diam(G)< 2^{-N}.$
It follows that given the set $\{B^1_1,\ldots,B^1_{Z_1}\}$, 
$m_{2N,N}^1,\ldots,m_{2N,N}^{Z_1}$ are all independent. 

Secondly, the random collection of boxes
$m_{2N,N}^i$ have the same distribution (after translation) as 
$m_{2N,N}^j$ for any $i,j\in\{1,\ldots,Z_1\}$ such that $i \neq j.$
Furthermore, $(m_{2N,N}^i)_{2^{N}}$ (recall the notation in \eqref{eqn:CDs}) 
has the same distribution as $m_N$
%$\{2^{-N}B:B\in m_N\}$ 
(again after translation).
Let $W_{2,i}$ be a maximal subset of $m_{2N,N}^i$ (using the same rule 
as the one used for picking $W_1$) such that 
$d(B^2_{j_1},B^2_{j_2})\geq 2^{-2N}$ for every 
$B^2_{j_1},B^2_{j_2}\in m^i_{2N,N}$ 
with $j_1 \neq j_2.$ Then, let $Z_{2,i}:=|W_{2,i}|$. 
It follows that the sequence $\{Z_{2,i}\}_{i=1}^{Z_1}$ is i.d.d.
and that every $Z_{2,i}$ has the same distribution as $Z_1$. 

Then, let $W_2 = \cup_{i=1}^{Z_1}W_{2,i}$, 
$Z_2=\sum_{i=1}^{Z_1}Z_{2,i}=|W_2|$ and
\[
\CW_2:=\bigcup_{B\in W_2} B.
\]
We observe that $Z_2 \leq \sum_{i=1}^{Z_1}|m_{2N,N}^i|\leq |m_{2N}|$,
and that $\CW_2 \subset \CC_{2N}$. Thus,  
$Z_2$ is the number of individuals in the second generation
of our branching process while the boxes of $W_2$ make up the second
generation of the corresponding family tree.

It is clear that we can proceed with the construction for further generations
in essentially the same way.
Of course, the general step in the construction is very similar to the steps 
above, 
so we shall be brief. Assume therefore that $(W_l,\CW_l,Z_l)$ has been constructed. 
First, we let $\{B^l_1,\ldots,B^l_{Z_l}\}$ be an 
enumeration of the level $lN$-boxes in $W_l$.
We think of $\{B^l_1,\ldots,B^l_{Z_l}\}$ as the offspring of generation $l.$  
Then, let 
\[
\CB_{(l+1)N}(B^l_i)=\{B \in \CB_{(l+1)N}: B \subset B^l_i\}
\]
which are the level $(l+1)N$-boxes that sits inside $B^l_i.$
Furthermore, we let
\[
m_{(l+1)N,lN}^i:=\{B\in \CB_{(l+1)N}(B^l_i):
\not \exists \ G \in \Phi_{(l+1)N}^{lN} \textrm{ such that } B\cap G\neq \emptyset \}.
\]

As before, the sequence $(m_{(l+1)N,lN}^i)_{1 \leq i \leq Z_l}$
consists of independent random variables. Then, we define
$W_{l+1,i}$ for $i\in \{1,\ldots,Z_l\}$ as above and let 
$Z_{l+1,i}:=|W_{l+1,i}|$. Finally, we let $W_{l+1}=\cup_{i=1}^{Z_l} W_{l,i}$,
$Z_{l+1}:=\sum_{i=1}^{Z_l}Z_{l+1,i}=|W_{l+1}|$ and 
\[
\CW_{l+1}:=\bigcup_{B\in W_l} B.
\]
Again, $\CW_{l+1} \subset \CC_{(l+1)N}.$

The sequence $(Z_l)_{l \geq 1}$ describes
a branching process with mean offspring distribution larger than $2^{\alpha N}$.
Since this is a supercritical process, it follows by standard theory that 
\[
\BP\left(\lim_{l\to \infty} Z_l=\infty\right)>0.
\]
On this event, we have that $\CW_{l} \neq \emptyset$ for every $l.$ 
Furthermore, the sets $\CW_l$ are all compact and $\CW_{l} \supset \CW_{l+1}$ 
for every $l\geq 1,$ and so
\[
\CC=\bigcap_{n=1}^\infty \CC_n=\bigcap_{l=1}^\infty \CC_{lN}
\supset \bigcap_{l=1}^\infty \CW_l \neq \emptyset.
\]
\fbox{}\\
{\bf Remark:} In \cite{BJT}, the corresponding result for the fractal ball model
was proved using a second moment method. That is, it was proved that for 
$\lambda<\frac{d \log 2}{\mu(\CA_1)}$ (see also Example \ref{ex:ball} in 
Section \ref{sec:examples})
there exists $c>0$ such that for every $n\geq 1,$
\[
\BP(|m_n| >0)\geq \frac{\BE[|m_n|]^2}{\BE[|m_n|^2]}\geq c>0.
\]
Using Fatou's lemma, this then implies that with positive probability, 
$|m_n|>0$ for every $n,$ and so $\CC_n \cap [0,1]^d\neq \emptyset$ 
for every $n.$ Furthermore, in that paper the balls were open, and so
$\CC_n\cap [0,1]^d$ was compact. Therefore, one can conclude that 
$\CC\cap [0,1]^d \neq \emptyset.$ 
This part of the argument will not work
in the more general setup of this paper. Indeed, if $\CC_n\cap [0,1]^d$ 
is not compact, the
conclusion that $\CC\cap [0,1]^d \neq \emptyset$ does not follow. 
Observe also that in general,
\[
\bigcup_{B \in m_{n_2}} B \not \subset  \bigcup_{B \in m_{n_1}} B
\]
whenever $n_1<n_2$. For instance, it is entirely possible that 
$m_{1}=\{[0,2^{-1}]^d\}$ while $m_2=\{[2^{-1},2^{-1}+2^{-2}]^d\}$.
Therefore, the fact that $|m_n|>0$ for every $n$ does not directly imply 
that $\CC\cap [0,1]^d \neq \emptyset$. Instead, we need to 
show the existence of nested, non-empty and closed subsets of 
$\CC_n \cap [0,1]^d$, which is 
why we need the branching process construction above. 

\bigskip
From the proof of Theorem \ref{thm:lbound} we conclude the following
corollary which may be of independent interest.
\begin{corollary} \label{corr:BranchingNumber}
For fixed $\lambda<\frac{d \log 2}{\mu(\CA_1)}$ and any $\epsilon>0,$ 
we can with positive probability 
embed a branching process with growth rate at least 
\[
2^{d-\lambda \frac{\mu(\CA_1)+\epsilon}{\log 2}},
\]
within $(m_n)_{n \geq 1}$.
\end{corollary}

\subsection{Upper bound}
The purpose of this section is to prove that 
$\lambda_e \leq \frac{d \log 2}{\mu(\CA_1)}$ which is then combined 
with Theorem \ref{thm:lbound} to prove Theorem \ref{thm:determinelambdae}.

\begin{theorem} \label{thm:ubound}
For any fractal model satisfying Assumptions \ref{assumption:reg1},
\ref{assumption:reg2} and \ref{assumption:reg3}
we have that 
\[
\lambda_e \leq \frac{d \log 2}{\mu(\CA_1)}.
\]
\end{theorem}
\noindent
{\bf Proof.}
Start by noting that an open ball of radius $2^{-n-1}$ can be inscribed in 
the box $[0,2^{-n}]^d$.
Therefore, 
\begin{eqnarray} \label{eqn:Mn1}
\lefteqn{\BP(B\in M_n)=\BP\left([0,2^{-n}]^d\in M_n\right)}\\
& & \leq \exp\left(-\lambda \mu(\{G \in \FG_n: D(o,2^{-n-1})\subset G\})\right). \nonumber
\end{eqnarray}
Then, let 
\[
S(G,r):=\{x\in G:d(x,\partial G)\geq r\},
\]
so that $S(G,r)$ is a shrunken version of $G.$ 
Obviously, $D(o,r) \subset G$ implies that $o \in G$ and 
$d(o, \partial G) \geq r$ so that $o \in S(G,r)$. It is also 
easy to see that the reverse implications are true. 
Thus,
\begin{eqnarray} \label{eqn:Mn2}
\lefteqn{\mu(\{G\in \FG_n: D(o,2^{-n-1})\subset G\})
=\mu(\{G\in \FG_n: o\in S(G,2^{-n-1})\})}\\
& & =\mu_p(\{G\in \FG_n: o\in S(G,2^{-n-1})\})
=\int_{\FG_n} I(o\in S(G,2^{-n-1})) d\mu_p(G), \nonumber
\end{eqnarray}
where the second equality follows since any $G\in \FG_n$ such that 
$\CL(G)=0$ must have an empty interior.
Next, similar to the proof of Lemma \ref{lemma:small_lambda}, we have that 
\begin{equation} \label{eqn:Mn3}
\int_{\FG_n}I(o\in S(G,2^{-n-1}))d\mu_p(G)
=\int_{\FG_n} \frac{\CL(S(G,2^{-n-1}))}{\CL(G)} I(o\in G)d\mu_p(G).
\end{equation}

Combining \eqref{eqn:Mn1}, \eqref{eqn:Mn2} and \eqref{eqn:Mn3} we then get that
\begin{eqnarray} \label{eqn:BMnineq}
\lefteqn{\BP(B\in M_n)
\leq \exp\left(-\lambda \int_{\FG_n} \frac{\CL(S(G,2^{-n-1}))}{\CL(G)} 
I(o\in G)d\mu_p(G)\right)}\\
& & =\exp\left(-\lambda \int_{\FG^o_n} \frac{\CL(S(G,2^{-n-1}))}{\CL(G)} 
d\mu_p(G)\right)
=\exp\left(-\lambda \sum_{l=1}^n \int_{\CA_l} \frac{\CL(S(G,2^{-n-1}))}{\CL(G)} 
d\mu_p(G)\right). \nonumber
\end{eqnarray}
Let 
\[
b_{l,n}:=\int_{\CA_l} \frac{\CL(S(G,2^{-n-1}))}{\CL(G)} d\mu_p(G)
\]
and observe that $b_{l,n}=b_{l-1,n-1}$ in the same way that we proved 
\eqref{eqn:alnrel}.
Thus, $b_{l,n}=\cdots=b_{1,n-l+1}$ for $l \leq n.$ Next, we have that 
\[
\bigcup_{n=1}^\infty S(G,2^{-n-1})=G \setminus \partial G.
\]
By the monotone convergence theorem and using Assumption \ref{assumption:reg2},
we then get that 
\begin{eqnarray*}
\lefteqn{\lim_{n \to \infty} b_{1,n}
=\int_{\CA_1} \lim_{n \to \infty}\frac{\CL(S(G,2^{-n-1}))}{\CL(G)} d\mu_p(G)}\\
& & =\int_{\CA_1} \frac{\CL(G \setminus \partial G)}{\CL(G)} d\mu_p(G)
=\int_{\CA_1} \frac{\CL(G \setminus \partial G)}{\CL(G)} I(\CL(\partial G)=0)
d\mu_p(G)
=\mu_p(\CA_1)
=\mu(\CA_1),
\end{eqnarray*}
where the last equality was shown in the proof of 
Lemma \ref{lemma:small_lambda}. We conclude that 
$b_{1,n}\uparrow \mu(\CA_1).$

Consider now some $\lambda>\frac{d \log 2}{\mu(\CA_1)}$ and let 
$\epsilon>0$ be such that also  
\begin{equation} \label{eqn:lambdalb}
\lambda>\frac{d \log 2}{\mu(\CA_1)-\epsilon}.
\end{equation}
There exists an $N_1<\infty$ such that 
$b_{1,n}\geq \mu(\CA_1)-\epsilon/2$ for every $n \geq N_1.$ Furthermore, 
there exists $N \geq N_1$ such that 
\[
(N-N_1)(\mu(\CA_1)-\epsilon/2) \geq 
N(\mu(\CA_1)-\epsilon).
\]
We then see that for every $n \geq N,$
\begin{eqnarray} \label{eqn:blnsum}
\lefteqn{\sum_{l=1}^n b_{l,n}=\sum_{l=1}^n b_{1,n-l+1}
=\sum_{l=1}^n b_{1,l}\geq \sum_{l=N_1+1}^{n} b_{1,l} }\\
& & 
\geq \sum_{l=N_1+1}^{n} (\mu(\CA_1)-\epsilon/2) 
=(n-N_1)(\mu(\CA_1)-\epsilon/2) 
\geq n(\mu(\CA_1)-\epsilon). \nonumber
\end{eqnarray}
Using \eqref{eqn:BMnineq} and \eqref{eqn:blnsum} we have that for 
every $n \geq N$, 
\begin{eqnarray*}
\lefteqn{\BE[|M_n|]
\leq 2^{dn}\exp\left(-\lambda \sum_{l=1}^n \int_{\CA_l} \frac{\CL(S(G,2^{-n-1}))}{\CL(G)} d\mu_p(G)\right)}\\
& & =\exp\left(nd\log 2-\lambda \sum_{l=1}^n b_{l,n}\right)
\leq \exp\left(nd\log 2-\lambda n(\mu(\CA_1)-\epsilon)\right)\\
& & =\exp\left(n (d\log 2-\lambda (\mu(\CA_1)-\epsilon))\right)
\longrightarrow 0,
\end{eqnarray*}
since $d \log 2-\lambda (\mu(\CA_1)-\epsilon)<0$ by \eqref{eqn:lambdalb}.

It follows that we must have that $\BP(\lim_{n \to \infty} |M_n|=\infty)=0,$
and so by Lemma \ref{lemma:survMninf} we can then conclude that 
\[
\BP(\CC(\lambda)\cap [0,1]^d =\emptyset)=1,
\]
from which the statement follows. 
\fbox{}\\

%We can now prove Theorem \ref{thm:determinelambdae}.\\
\noindent
{\bf Proof of Theorem \ref{thm:determinelambdae}.}
This is immediate from Theorems \ref{thm:lbound} and \ref{thm:ubound}.
\fbox{}\\

\section{Proof of Theorem \ref{thm:dieoutatcrit}} \label{sec:proofcrit}
The purpose of this section is to prove Theorem \ref{thm:dieoutatcrit}
and also to discuss the connection between condition \eqref{eqn:extracond} 
and the notion of thin as introduced in \cite{NW}. In addition, we shall 
derive an alternative expression for $\lambda_e.$ %and $\dim_H(\CC(\lambda))$.

\medskip

We start by proving Theorem \ref{thm:dieoutatcrit} which basically just combines 
the proof of Theorem \ref{thm:ubound} with condition 
\eqref{eqn:extracond}. \\
\noindent
{\bf Proof of Theorem \ref{thm:dieoutatcrit}.}
Note that by the definition of $[\partial G]^r,$ and since 
$\mu_p(\CA_1)=\mu(\CA_1),$
\begin{eqnarray*}
\lefteqn{\mu(\CA_1)-b_{1,l}=\int_{\CA_1}1-\frac{\CL(S(G,2^{-l-1}))}{\CL(G)} 
d\mu_p(G)}\\
& & =\int_{\CA_1}\frac{\CL(G \setminus S(G,2^{-l-1}))}{\CL(G)} d\mu_p(G)
=\int_{\CA_1}\frac{\CL([\partial G]^{2^{-l-1}})}{\CL(G)} d\mu_p(G).
\end{eqnarray*}
Furthermore, 
\begin{eqnarray*}
\lefteqn{\sum_{l=1}^\infty \int_{\CA_1}
\frac{\CL([\partial G]^{2^{-l-1}})}{\CL(G)} d\mu_p(G)
\leq \sum_{l=1}^\infty \int_{\CA_1}
\frac{\CL([\partial G]^{2^{-l}})}{\CL(G)} d\mu_p(G)}\\
& & = \int_{\CA_1} \sum_{l=1}^\infty 2\int_{2^{-l}}^{2^{-l+1}} 
\frac{1}{2^{-l+1}}
\frac{\CL([\partial G]^{2^{-l}})}{\CL(G)} dr d\mu_p(G)
\leq 2 \int_{\CA_1} \int_0^1 
\frac{1}{r} \frac{\CL([\partial G]^r)}{\CL(G)} dr d\mu_p(G)<\infty,
\end{eqnarray*}
where the last inequality is condition \eqref{eqn:extracond}.
Therefore, as in the proof of Theorem \ref{thm:ubound}, we have that 
for $\lambda=\lambda_e$
\begin{eqnarray} \label{eqn:EMncrit}
\lefteqn{\BE[|M_n|]
\leq \exp\left(d n\log 2 -\lambda_e\sum_{l=1}^n b_{1,l}\right)}\\
& &=\exp\left(\frac{d \log 2}{\mu(\CA_1)}
\left(n\mu(\CA_1)-\sum_{l=1}^n b_{1,l}\right)\right) \nonumber \\
& & =\exp\left( \frac{d \log 2}{\mu(\CA_1)}\sum_{l=1}^n \int_{\CA_1}
\frac{\CL([\partial G]^{2^{-l-1}})}{\CL(G)} d\mu_p(G)\right) \nonumber  \\
& & \leq \exp\left( \frac{d \log 2}{\mu(\CA_1)}\sum_{l=1}^\infty \int_{\CA_1}
\frac{\CL([\partial G]^{2^{-l-1}})}{\CL(G)} d\mu_p(G)\right)<\infty. \nonumber
\end{eqnarray}
We conclude that 
$\limsup_{n \to \infty}\BE[|M_n|]<\infty$, and the statement follows by using 
Lemma \ref{lemma:survMninf} as before.
\fbox{}\\
{\bf Remark:} 
It is clear from the proofs of Theorems \ref{thm:determinelambdae} 
and \ref{thm:dieoutatcrit}, that for a fractal model $\mu$ 
satisfying Assumptions \ref{assumption:reg1}, \ref{assumption:reg2} 
and \ref{assumption:reg3}, $\mu_0$ is irrelevant when it comes 
to the question of covering subsets of $\BR^d$. However, this is 
not apriori obvious, which is why we chose not to restrict our
attention to measures $\mu$ such that $\mu=\mu_p$.

\medskip

In order to proceed, we will consider a particular way of constructing 
scale invariant measures $\mu.$ 
We shall be somewhat informal, the missing details can be found in \cite{NW}. 
We note that the models considered in \cite{NW} do not have a cutoff as we 
have imposed on $\mu$ here. Instead, they consider the model in a bounded domain 
$D$ where any set $G\in \Phi_\lambda(\mu)$ such that 
$G \cap (\partial D \cup D^c) \neq \emptyset$ 
is removed. However, this has little impact on the construction.
Let $L(G)$ denote the center of mass of $G$. 
%Firstly, for any closed set $G$ we let 
%$L_1(G):=\min \{x_1:(x_1,\ldots,x_d)\in G\},$ and inductively we then let 
%\[
%L_k(G):=\min \{x_k:(L_1(G),\ldots,L_{k-1}(G),x_k,\ldots,x_d)\in G\}
%\] 
%for 
%$k \leq d.$ Then, we define $L(G):=(L_1(G),\ldots,L_d(G))$. 
This creates 
a sort of ''anchor'' for the set $G$, and in fact any such anchor would do.

Let $\nu_o$ be a finite measure on sets $H\subset \BR^d$ such that 
$\diam(H)=1$ and $L(H)=o$, and consider the product measure 
\[
\nu_o \times \rho^{-(d+1)}d\rho \times \CL(dx).
\]
Using the transformation $G=\rho H+x$ we can define a scale invariant measure 
$\mu$ with cutoff 1, through 
\[
\int F(G) d\mu(G)
:=\int_{\BR^d}\int_0^1 \int F(\rho H+x)d\nu_o\times \rho^{-(d+1)}d\rho\times \CL(dx),
\]
where $F$ is in some suitable class of functions.
Furthermore, any scale invariant measure $\mu$ can be constructed in such a way 
(see again \cite{NW}). There is a tiny difference in that we assume that $\mu$
is a measure on $\{G \subset \BR^d: \diam (G)<1\}$ rather than 
$\{G \subset \BR^d: \diam (G)\leq 1\}$. However, this is only a technical 
nuisance. The condition of $\mu$ being thin is defined in \cite{NW} to mean that 
\begin{equation} \label{eqn:thin}
\int_0^1 \frac{1}{r} \int \CL(E(H,r)\setminus S(H,r))d\nu_o(H) dr<\infty. 
\end{equation}
%It is immediate from \eqref{eqn:thin} that $\mu$ being thin implies 
%Assumption \ref{assumption:reg2}. 

%WE HAVE THAT $\nu_o$ being finite implies that 
We now turn to the proof of Proposition \ref{prop:condassrel} which relates
Assumption \ref{assumption:reg2} with conditions \eqref{eqn:extracond}
and \eqref{eqn:thin}.

\noindent
{\bf Proof of Proposition \ref{prop:condassrel}.} 
We start by proving $(i).$
If $\mu$ is thin, it follows from \eqref{eqn:thin} that 
$\int\CL(\partial H) d\nu_o(H)=0.$ Therefore $\nu_o$ is supported on sets 
such that $\CL(\partial H)=0$. It follows that this is also the case for 
the measure $\mu$, and so Assumption \ref{assumption:reg2} holds.

In order to show that $\mu$ being thin 
implies condition \eqref{eqn:extracond} we observe that
\begin{eqnarray} \label{eqn:relcond1}
\lefteqn{\int_{\CA_1}\frac{\CL([\partial G]^r)}{\CL(G)}d\mu_p(G)
=\int_{\CA_1}\frac{\CL([\partial G]^r)}{\CL(G)}I(\CL(G)>0)d\mu(G)}\\
%=\int \frac{\CL(G\setminus S(G,2^{-l}))}{\CL(G)} I(o \in G)d\mu(G)}\\
& & =\int_{1/2}^1\int_{\BR^d} \int 
\frac{\CL((\rho H+x) \setminus S(\rho H+x,r))}{\CL(\rho H+x)}  
I(o\in \rho H+x) I(\CL(H)>0)\frac{1}{\rho^{d+1}}d\nu_o(H)\CL(dx)d\rho\nonumber\\
& & =\int_{1/2}^1\int_{\BR^d} \int 
\frac{\CL(\rho H\setminus S(\rho H,r))}{\CL(\rho H)}  
I(o\in \rho H+x) I(\CL(H)>0)
\frac{1}{\rho^{d+1}}d\nu_o(H) \CL(dx)d\rho \nonumber \\
& & =\int_{1/2}^1\int 
\frac{\CL(H\setminus S(H,r/\rho))}{\CL(H)}\CL(\rho H) I(\CL(H)>0)
 \frac{1}{\rho^{d+1}}d\nu_o(H) d\rho \nonumber  \\
&& =\int_{1/2}^1\int \CL([\partial H]^{r/\rho})
 \frac{1}{\rho}d\nu_o(H) d\rho 
 \leq \log 2 \int \CL([\partial H]^{2r})d\nu_o(H). \nonumber
\end{eqnarray}
Thus, 
\begin{eqnarray*}
\lefteqn{\int_{\CA_1}\int_0^1\frac{1}{r}\frac{\CL([\partial G]^r)}{\CL(G)} 
dr d\mu_p(G)
\leq \log 2 \int_0^1\frac{1}{r} \int \CL([\partial H]^{2r})d\nu_o(H) dr} \\
& & \leq C+\log 2 \int_0^{1/2}\frac{1}{r} \int \CL([\partial H]^{2r})d\nu_o(H) dr \\
& & = C+\log 2 \int_0^1 \frac{1}{r} \int \CL([\partial H]^{r})d\nu_o(H)dr \\
& & \leq C+\log 2 \int_0^1 \frac{1}{r} \int \CL(E(H,r)\setminus S(H,r))d\nu_o(H)dr,
\end{eqnarray*}
and so we see that \eqref{eqn:thin} indeed implies condition \eqref{eqn:extracond}.

We now turn to $(ii)$. By a slight adjustment we can reverse the
inequality in \eqref{eqn:relcond1}. Furthermore, if the measure $\mu$ is supported 
on compact sets, then $[\partial H]^{r} \supset \partial H$ for every $r\geq 0,$
and so
\[
\int_{\CA_1} \frac{\CL([\partial G]^r)}{\CL(G)} d\mu_p(G)
\geq \log 2\int \CL([\partial H]^{r})d\nu_o(H)
\geq \log 2\int \CL(\partial H)d\nu_o(H).
\]
Hence, 
\[
\log 2 \int_0^1\frac{1}{r} \int \CL(\partial H)d\nu_o(H)dr
\leq \int_0^1\frac{1}{r} \int_{\CA_1} \frac{\CL([\partial G]^r)}{\CL(G)} d\mu_p(G)dr,
\]
and so if condition \eqref{eqn:extracond} holds, we must have 
that $\int \CL(\partial H)d\nu_o(H)=0$, which implies Assumption 
\ref{assumption:reg2}.

In order to prove that condition \eqref{eqn:extracond} does not imply 
Assumption \ref{assumption:reg2} in general, we will construct a counterexample.
%Let therefore $\nu_o$ be concentrated on the following set.
Let $x_1,x_2,\ldots$ be an enumeration of the rational points in $[0,1]^d,$
where $d\geq 2,$
and let ${\rm Box}(x,r)$ denote an open box centred at $x$ with side length $r.$ 
Then, we define 
\[
H:=\bigcup_{k=1}^\infty {\rm Box}(x_k,2^{-k-1}).
\]
We see that $H$ is open and that $\CL(H)\leq 2^{-d}$. 
We pick our measure $\nu_o$ as being unit mass on the set $H.$ 
Obviously $\diam(H)>1,$ but this is easily rectified and left for the reader.
We have that $\partial H \supset [0,1]^d \setminus H$ so that 
$\CL(\partial H)\geq 1-2^{-d}$ and so Assumption \ref{assumption:reg2} cannot 
hold for this choice of $H$.

Observe that 
\[
\CL([\partial {\rm Box}(x,s)]^r)=\left\{\begin{array}{ll}
s^d & \textrm{if } s\leq 2r \\
s^d-(s-2r)^d & \textrm{if } s>2r.
\end{array} \right.
\]
Observe also that since the boxes ${\rm Box}(x_k,2^{-k-1})$ and
${\rm Box}(x_l,2^{-l-1})$ can overlap for $k \neq l$, there exists 
some $c=c(d)<\infty$ independent of $r$ and such that 
\begin{eqnarray*}
\lefteqn{\CL([\partial H]^{r})
\leq \sum_{k=1}^\infty \CL([\partial {\rm Box}(x_k,2^{-k-1})]^{r})}\\
& & \leq \sum_{k=1}^{\lceil -\log_2 r \rceil-3} 2^{-dk-d}-(2^{-k-1}-2r)^d
+\sum_{k=\lceil -\log_2 r \rceil-2}^\infty 2^{-dk-d}\\
& & \leq \frac{cr}{2}
+\sum_{k=1}^{\lceil -\log_2 r \rceil-3} 2^{-dk-d}(1-(1-r2^{k+2})^d)
\leq \frac{cr}{2}+\sum_{k=1}^{\lceil -\log_2 r \rceil-3} 2^{-dk-d} d r 2^{k+2}
\leq cr,
\end{eqnarray*}
where we use that $d\geq 2$. Furthermore, we also used 
that if $k \leq \lceil -\log_2 r \rceil-3$ 
then $r 2^{k+2}\leq r 2^{\lceil -\log_2 r \rceil-1}\leq 1$ 
and that $1-(1-x)^d \leq dx$ whenever $0\leq x \leq 1.$
Obviously, we here have that $\mu=\mu_p$ and
using \eqref{eqn:relcond1} we see that  
\[
\int_0^1 \frac{1}{r} \int_{\CA_1} 
\frac{\CL([\partial G]^r)}{\CL(G)}d\mu_p(G) dr
\leq \log 2 \int_0^1 \frac{1}{r} \int \CL([\partial H]^{2r})d\nu_o(H)
<\infty,
\]
and so condition \eqref{eqn:extracond} holds.

Finally, we prove $(iii)$ by constructing an example satisfying Assumption 
\ref{assumption:reg2} but not condition \eqref{eqn:extracond}. This construction 
will be done sequentially, i.e. we will construct closed sets 
$K_1,K_2,\ldots \subset [0,1]^d$
and then define $H:=\left(\cap_{n=1}^\infty K_n\right)^c.$ Every set $K_n$
will be the union of a set of boxes $\CK_n \subset \CB_n$ so that 
by definition of $\CB_n,$
\[
K_n:=\bigcup_{B\in \CK_n}B \subset [0,1]^d.
\]
Let $\CK_1=\CB_1$ so that $K_1=[0,1]^d$. Then, let $\CK_2$ be the set of boxes
$B\in \CB_2$ such that $B=x+[0,2^{-2}]^d$ where $x_1+\cdots+x_d$ is an even 
number. Note that $\CK_2$ consists of $2^{2d}/2$ boxes in a 
($d$-dimensional) chequerboard pattern.
In general, let $\CK_{n+1}$ consist of $\lceil 2^{d(n+1)}/(n+1) \rceil$ boxes 
of $\CB_{n+1}$, picked so that $K_{n+1}\subset K_n,$ and such that 
for any $B\in \CK_{n},$ there exists $B'\in \CK_{n+1}$
such that $B' \subset B$. This condition makes sure that the box 
$B\in \CK_{n}$ is not removed in its entirety. 
It is straightforward to see that this is possible 
as we now (somewhat informally) explain. Observe that the set $K_n$ can be written 
as a union of $2^d \lceil 2^{dn}/n \rceil$ boxes from $\CB_{n+1}$. 
Therefore, the number of boxes that needs to be deleted in order 
to create $K_{n+1}$ can be bounded by 
\[
2^d\lceil 2^{dn}/n \rceil-\lceil 2^{d(n+1)}/(n+1) \rceil
\leq 2^{d(n+1)}\left(\frac{1}{n}-\frac{1}{n+1}\right)+1
=\frac{2^{d(n+1)}}{n(n+1)}+1\leq \frac{1}{2}2^d \lceil 2^{dn}/n \rceil,
\]
where the last inequality holds for $n\geq 2$ and $d\geq 2.$ However, 
a similar calculation for $n=1$ shows that for any $d\geq 2$ and any $n\geq 1,$
\[
2^d\lceil 2^{dn}/n \rceil-\lceil 2^{d(n+1)}/(n+1) \rceil
\leq \frac{1}{2}2^d \lceil 2^{dn}/n \rceil.
\]
Therefore, in every step of the construction, 
at least half of the boxes are kept and so it is always possible to keep 
at least one level $(n+1)$ subbox of $B\in \CK_n.$

We then note that since $H^c=\bigcap_{n=1}^\infty K_n$ is closed, 
we must have that 
\[
\CL(\partial H)=\CL(H^c)=\lim_{n \to \infty}\CL(K_n)
=\lim_{n \to \infty}\lceil 2^{dn}/n \rceil 2^{-dn}=0,
\]
and so Assumption \ref{assumption:reg2} is satisfied. We also note that 
for every $B \in \CK_n$ there exists $x\in H^c$ such that $x\in B$. 
This follows from the constraint in the construction of $\CK_{n+1}$ from
$\CK_{n}$ which implies that $B \cap K_n \neq \emptyset$ for every 
$n \geq 1,$ and so $B \cap H^c \neq \emptyset$ by compactness. This in turn 
implies that $d(y,H^c)\leq \sqrt{d}2^{-n}$ for every $y \in K_n$
and so $[\partial H]^{\sqrt{d}2^{-n}} \supset K_n\setminus H^c$ 
for every $n\geq 1.$  Thus, for $\sqrt{d}2^{-n} \leq r \leq \sqrt{d}2^{-n+1}$
\[
\CL([\partial H]^r)
\geq \CL([\partial H]^{\sqrt{d}2^{-n}}) \geq \CL(K_n \setminus H^c)
=\CL(K_n)=\lceil 2^{dn}/n \rceil 2^{-dn} \geq 1/n.
\]
Therefore, by letting $\nu_o$ be unit mass on $H$ we see that 
\[
\int_0^1 \frac{1}{r}\int \CL([\partial H]^r) d\nu_o(H)dr
\geq \sum_{n=d}^\infty \int_{\sqrt{d}2^{-n}}^{\sqrt{d}2^{-n+1}}
\frac{1}{\sqrt{d}2^{-n+1}}\frac{1}{n}dr
=\sum_{n=d}^\infty\frac{1}{2n}=\infty,
\]
and so condition \eqref{eqn:extracond} is not satisfied.
\fbox{}\\

\noindent
{\bf Remarks:} Informally, 
the example in part $(ii)$ of Proposition \ref{prop:condassrel} works 
since the the set $H$
has a small inner boundary but a large outer boundary.

Similar to \eqref{eqn:relcond1}, one obtains that 
\begin{eqnarray*}
\lefteqn{\mu(G:G \cap D(o,1)\neq\emptyset,\delta\leq \diam(G)\leq 1)}\\
& & =\int_\delta^1 \int_{\BR^d} \int I(o \in E(\rho H+x,1))
\frac{1}{\rho^{d+1}}d\nu_o(H)\CL(dx)d\rho \\
& & =\int_\delta^1\int \CL(E(H,1/\rho))\frac{1}{\rho}d\nu_o(H)d\rho,
\end{eqnarray*}
which is finite whenever $\nu_o$ is finite. Thus, the finiteness of $\nu_o$
is equivalent to Assumption \ref{assumption:reg1}.
However, if we allow ourselves to consider 
also infinite measures $\nu_o,$ it is possible to find an example so that 
condition \eqref{eqn:extracond} is satisfied while Assumption \ref{assumption:reg1} 
is not. We will give an informal description of how this is done.
Let $\nu_o$ be a measure on sets $[0,1]\times [0,h]$ 
%using $h^{-3/2}dh$ as a measure on $(0,1]$ so that 
given by $\nu_o(\CL(H)\geq l)=\int_l^1 h^{-3/2}dh$.
Obviously, $\nu_o$ is then an infinite measure. However, we see that
\[
\int \CL([\partial H]^{r})d\nu_o(H) 
\leq \int_0^{2r}h^{-1/2}dh+4\int_{2r}^1 r h^{-3/2}dh 
\leq C \sqrt{r},
\] 
%\[
%\int \CL(H\setminus S(H,2^{-l}))d\nu_o(H) 
%\leq \int_0^{2^{-l}}h^{-1/2}dh+4\int_{2^{-l}}^1 2^{-l}h^{-3/2}dh 
%\leq C 2^{-l/2},
%\] 
and as above we get that condition \eqref{eqn:extracond} is satisfied.

%\begin{eqnarray*}
%\lefteqn{\int \CL(H\setminus S(H,2^{-l}))d\nu_o(H) }\\
%& & \leq \int_0^{2^{-l}}h^{-1/2}dh+4\int_{2^{-l}}^1 2^{-l}h^{-3/2}dh 
%\leq C 2^{-l/2}
%\end{eqnarray*}

\bigskip

We now proceed to find the canonical expression for $\lambda_e$ mentioned 
in the Introduction. To that end, observe that as above, 
\begin{eqnarray} \label{eqn:muA1alt}
\lefteqn{\mu(\CA_1)
=\int_{1/2}^1\int_{\BR^d}\int I(o\in \rho H+x)\frac{1}{\rho^{d+1}}
d\nu_o(H)\CL(dx)d\rho}\\
& & =\int_{1/2}^1 \int \CL(H) \frac{1}{\rho} d\nu_o(H)d\rho
=\log 2\int \CL(H) d\nu_o(H)
=\nu_o(\CL(H))\log 2. \nonumber
\end{eqnarray}
Thus we may write $\lambda_e$ as
\begin{equation} \label{eqn:lambdae_alt}
\lambda_e=\frac{d}{\nu_o(\CL(H))}.
\end{equation}
This expression relates the critical value $\lambda_e$ to the $\nu_o$-measure of 
the volume of sets of diameter one. Of course, while $\nu_o$ must be finite by 
Assumption \ref{assumption:reg1}, it does not have to be a probability measure.
%Therefore, $\nu_o(\CL(H))$ is not quite the expected volume of the set
%$H$ but rather a quasi-expectation.

We end this section by giving the proof of Corollary \ref{corr:closedopensame}. \\
\noindent
{\bf Proof of Corollary \ref{corr:closedopensame}.}
It follows from Assumption \ref{assumption:reg2} that
$\nu_o(\CL(H))=\nu_o(\CL(H \cup \partial H))$ and so 
the statement follows immediately from \eqref{eqn:lambdae_alt}.
\fbox{}\\

\noindent
{\bf Remark:}
In a similar manner, if $\nu_o,\nu'_o$ are both point masses on two different 
sets having the same volume, then $\lambda_e(\nu_o)=\lambda_e(\nu'_o)$.

\section{Examples} \label{sec:examples}
In this section we will briefly review some examples. We will avoid delving 
into lengthy calculations and discussions and will not give precise references 
to well known results (e.g. the fractal dimension of the boundary of a von Koch
snowflake). Furthermore, although one could in principle just normalize the 
measure $\mu$ in every example (so that say $\mu(\CA_1)=1$) this would obviously
defeat the purpose. Therefore, the normalizations in the examples are what we 
consider to be canonical.

It is convenient to have the following proposition which makes it easier to 
verify condition \eqref{eqn:extracond}. Here, $\overline{\dim}_B(\partial H)$ 
denotes the so called upper box (or Minkowski) dimension of the set $\partial H$. 
We will not review this concept here but rather refer to Chapter 3
of \cite{Falconer}.
\begin{proposition} \label{prop:boxdim}
Assume that $\nu_o$ is concentrated on a finite collection 
of sets $H_1,\ldots,H_K$ such that 
$\max_{k=1,\ldots,K} \overline{\dim}_B(\partial H_k)<d$.
 Then the resulting measure $\mu$ is thin. 
\end{proposition}
\noindent
{\bf Proof.}
As before, we let $D(x,r)$ denote a ball of radius $r$ centred at $x.$ 
Let $\CP_r(H_k)$ denote a minimal collection of points such that 
\[
\partial H_k \subset \bigcup_{x \in \CP_r(H_k)} D(x,r).
\]
Since any point $y\in E(H_k,r)\setminus S(H_k,r)$ must be within distance 
at most $r$ from the boundary, it follows that $y\in D(x,2r)$ for some 
$x\in \CP_r(H_k).$ Thus,
\[
E(H_k,r)\setminus S(H_k,r) \subset \bigcup_{x \in \CP_r(H_k)} D(x,2r),
\]
and so $\CL(E(H_k,r)\setminus S(H_k,r))\leq |\CP_r(H_k)| (2r)^d$.
Let $\alpha$ be such that $\max_{k=1,\ldots,K} \overline{\dim_B(\partial H_k)}
<\alpha<d$. It follows that for any $r$ small enough we must have that 
$|\CP_r(H_k)|\leq r^{-\alpha}$ for $k=1,\ldots,K.$
Therefore, 
\[
\CL(E(H_k,r)\setminus S(H_k,r))\leq r^{-\alpha}(2r)^d=2^d r^{d-\alpha},
\]
and so we see that \eqref{eqn:thin} is satisfied.
\fbox{}\\

Our first example connects the results of this paper with previous research.
\begin{example} \label{ex:ball}
The fractal ball model.
\end{example}
\noindent
In this example we shall briefly consider the classical fractal Poisson ball model.
%The example shows that the results of \cite{BJT} and \cite{BE} are special cases
%of Theorems \ref{thm:determinelambdae} and \ref{thm:dieoutatcrit}.
Let $\nu(dr)=r^{-d-1}dr$ be a measure on $(0,1/2]$ and let $\mu=\CL \times \nu$
be the product measure on $\BR^d \times (0,1/2]$. We can view $\mu$ as a measure
on $\FG$ by associating a point $(x,r)\in \BR^d \times (0,1/2]$ with the ball
(open or closed) $D(x,r)$ of radius $r.$ 
The measure $\mu$ is then taken to be the intensity 
measure of the fractal ball model, and it is not hard to prove that it is 
indeed scale invariant and satisfies Assumptions \ref{assumption:reg1},
\ref{assumption:reg2} and \ref{assumption:reg3}. 
Furthermore, 
\[
\mu(\CA_1)=\int_{\BR^d} \int_{1/4}^{1/2} I(o\in D(x,r))\nu(dr)\CL(dx)
=\int_{1/4}^{1/2} v_d r^d r^{-d-1}dr=v_d\log 2,
\]
where $v_d$ denotes the volume of the unit ball in $\BR^d.$ Thus,
by Theorem \ref{thm:determinelambdae} we have that
\[
\lambda_e=\frac{d \log 2}{\mu(\CA_1)}=\frac{d}{v_d}.
\]
Alternatively, we can use the construction of Section \ref{sec:proofcrit}
with $\nu_o$ being unit mass on the unit ball to obtain the same result.

Using Proposition \ref{prop:boxdim} we see that $\mu$ is thin, and so 
by Proposition \ref{prop:condassrel}, condition \eqref{eqn:extracond} is 
satisfied. Thus, it follows from Theorem \ref{thm:dieoutatcrit} 
that the Poisson ball model satisfies $\CC(\lambda_e)=\emptyset$. This 
was previously proven in \cite{BJT} for the case of open balls.

\bigskip

Our next example considers a somewhat more irregular ''fundamental''
shape than balls. In particular, the shape is non-convex and has a fractal boundary.
\begin{example} \label{ex:snow}
The scale invariant von Koch snowflake model.
\end{example}
\noindent
An early example of a set with a fractal boundary is the so-called 
von Koch snowflake (see for instance \cite{Falconer} p. xiv  for a construction) 
which is a subset of
$\BR^2$. We will let $S(x,l,\theta)$ denote such a snowflake with diameter $l$
and orientation $\theta$ compared to some chosen base direction.
The boundary of the snowflake is known to have box dimension $\log 4/\log 3$ while 
the area of $S(x,l)$ equals $\frac{3\sqrt{3}}{10}l^2.$ We can easily turn 
this into a scale invariant snowflake model by a construction similar to 
Example \ref{ex:ball} as we now explain. Let $\Theta$ be uniform measure
on $[0,2\pi)$ and let $\nu(dl)=l^{-3}dl$ on $(0,1]$ (note that $l$ denotes
the diameter while $r$ in the previous example was the radius). Then, we let 
$\mu=\CL \times \nu \times \Theta$ be product measure on 
$\BR^d \times (0,1] \times [0,2\pi)$ and to any point 
$(x,l,\theta) \in \BR^d \times (0,1] \times [0,2\pi)$ we associate a snowflake
$S(x,l,\theta)$. Then, using $\mu$ as our intensity measure we have constructed
a scale invariant Poisson snowflake model. It is easy to check that Assumptions 
\ref{assumption:reg1}, \ref{assumption:reg2} and \ref{assumption:reg3} 
are all satisfied for this model. 

Furthermore we have that 
\begin{eqnarray*}
\lefteqn{\mu(\CA_1)=\int_{1/2}^1 \int_{\BR^2} \int_0^{2\pi} I(o \in S(x,l,\theta))
d\Theta \CL(dx) \nu(dl)}\\
& & =\int_{1/2}^1 \int_{\BR^2}I(o \in S(x,l,0))d\CL(x) \nu(dl)
=\int_{1/2}^1 \frac{3 \sqrt{3}}{10}l^2 l^{-3} dl
=\frac{3 \sqrt{3}}{10}\log 2,
\end{eqnarray*}
and so according to Theorem \ref{thm:determinelambdae} we obtain
\[
\lambda_e=\frac{2 \log 2}{\mu(\CA_1)}=\frac{20}{3 \sqrt{3}}.
\]
Here, we could also have used $\nu_o$ being a rotational invariant measure on 
snowflakes of diameter 1 centred at $o.$

In addition we can use Propositions \ref{prop:boxdim} (the fact that we add 
rotation to the snowflakes does not change anything) and \ref{prop:condassrel} 
to conclude that $\mu$ satisfies condition \eqref{eqn:extracond}. Then, 
Theorem \ref{thm:dieoutatcrit} tells us that for 
$\lambda=\lambda_e$, $\BR^2$ is almost surely covered by 
scale invariant snowfall.

\begin{example} \label{ex:BLS}
The Brownian loop soup
\end{example}
\noindent
The Brownian loop soup was introduced in \cite{LW} and has since received 
much attention as it is intimately connected to the theory of so-called 
Schramm-Loewner Evolution ($SLE$) and Conformal Loop Ensembles ($CLE$)
(see for instance \cite{Lawler}).
The Brownian loop measure $\mu^{loop}$ is a scale invariant measure on 
Brownian loops in the plane. In this example we will impose a cutoff on the 
diameter of the loops so that no loop has diameter larger than 1. 
Furthermore, given a loop $\gamma$ of the loop soup, we define the set $G$ 
as the complement of the unbounded component created by $\gamma$. 
The measure $\mu$ is then defined to be the measure on the sets $G$ 
induced by the outer boundaries of the loops as we just described.

It is known (see for instance \cite{Camia} Lemma B.8 p. 266) that 
in this case we have that 
\[
\mu(\CA_1)=\frac{1}{5}\log 2,
\]
from which it follows that 
\[
\lambda_e=\frac{2 \log 2}{\mu(\CA_1)}=10.
\]
Furthermore, it is proven in \cite{NW} Lemma 4 that while the Brownian loop
soup in itself is not thin, the measure $\mu$ on the induced sets $G$ is.
Therefore, also this example satisfies the conditions of
Theorem \ref{thm:dieoutatcrit} and so $\CC(\lambda_e)=\emptyset.$

\medskip

\noindent
{\bf Acknowledgements:} The author would like to thank F. Camia, J. Tykesson
and J. Steif for useful comments on an early version of this paper. In 
addition, the author is grateful towards the anonymous referee for many 
useful comments and suggestions.


\begin{thebibliography}{99}

\bibitem{BE} Bierm\'{e} H. and Estrade A. Covering the whole space with
Poisson random balls. {\em ALEA Lat. Am. J. Probab. Math. Stat.}
{\bf 9}, (2012), 213--229.


\bibitem{BC} Broman E. and Camia F.
Universal behavior of connectivity properties in fractal percolation models.
{\em Electron. J. Probab.} {\bf 15}, (2010), 1394--1414.


\bibitem{BJT} Broman E., Jonasson J. and Tykesson J.
The existence phase transition for two Poisson random fractal models. 
{\em Electron. Commun. Probab.} {\bf 22}, (2017), Paper No. 21.



\bibitem{Camia} Camia F. Scaling Limits, Brownian Loops and
Conformal Fields. {\em Advances in Disordered Systems, Random Processes and 
Some Application} Cambridge University Press (2017), 205--269. 
%Lemma B.8 is on page 266.



\bibitem{DM} Dekking, F. M. and Meester, R. W. J.
On the structure of Mandelbrot's percolation process and other random Cantor sets.
{\em J. Statist. Phys.} {\bf 58}, (1990), no. 5-6, 1109--1126.


\bibitem{Dvoretsky} Dvoretsky A. 
On covering a circle with randomly placed arcs
{\em Proc. Nat. Acad. Sci. U. S. A.}
{\bf 42}, (1956), 199--203.


\bibitem{Falconer} Falconer K. {\em Fractal Geometry}, John Wiley \& Sons, Ltd.,
(2014).

%\bibitem{GF} Garban C. and Ferreras J. The expected area of the filled Brownian 
%loop is $\pi/5.$

\bibitem{Lawler} Lawler G.F. {\em Conformally invariant processes in the plane},
Mathematical Surveys and Monographs, {\bf 114}, American Mathematical Society
(2005).

\bibitem{LW} Lawler G.F. and Werner W. The Brownian loop soup.
{\em Probab. Theory Relat. Fields}, {\bf 128}, (2004), 565--588.

\bibitem{LP} Lyons R. and Peres Y. {\em Probability on trees and networks}
Cambridge Series in Statistical and Probabilistic Mathematics, {\bf 42} (2016). 

\bibitem{Molchanov} Molchanov I. {\em Theory of Random Sets} 2nd Edition,
Springer-Verlag, (2017).

\bibitem{NW} Nacu S. and Werner W. Random soups, carpets and 
fractal dimensions. {\em J. Lond. Math. Soc.} (2) {\bf 83} 
(2011), no. 3, 789--809. 

\bibitem{Shepp} Shepp L.A. Covering the line with random intervals
{\em Z. Wahrscheinlichkeitstheorie und Verw. Gebiete}, {\bf 23}, (1972), 163--170.

%\bibitem{Werner} SLEs as boundaries of clusters of Brownian loops.

\end{thebibliography}
\end{document}